\documentclass[a4, 12pt, reqno]{svjour3}

\smartqed
\usepackage{pifont,geometry,latexsym,graphicx,mathptmx}
\usepackage{amssymb,amsmath,amsfonts,amscd,stmaryrd}
\usepackage{a4wide}

\newtheorem{Theorem}{Theorem} 

\newtheorem{Remark}[Theorem]{Remark}
\newtheorem{Lemma}[Theorem]{Lemma}

\newtheorem{Definition}[Theorem]{Definition}

\numberwithin{Theorem}{section}
\numberwithin{equation}{section}

\newcommand{\hh}{\hspace*{.48in}}
\newcommand{\h}{\hspace*{.24in}}
\newcommand{\hhh}{\hspace*{.72in}}

\begin{document}
\title{Adjoint methods for static Hamilton-Jacobi equations}
\author{Hung Vinh Tran*}
\institute{Department of Mathematics\\
University of California, Berkeley, CA 94720.\\
\email{tvhung@math.berkeley.edu}\\
*Supported in part by VEF fellowship.}
\date{} 
\maketitle
\begin{abstract}
We use the adjoint methods to study the static Hamilton-Jacobi equations and to prove the speed of convergence for those equations. The main new ideas are to introduce adjoint equations corresponding to the formal linearizations of regularized equations of vanishing viscosity type, and from the solutions $\sigma^{\epsilon}$ of those we can get the properties of the solutions $u$ of the Hamilton-Jacobi equations. We classify the static equations into two types and present two new ways to deal with each type. The methods can be applied to various static problems and point out the new ways to look at those PDE.
\keywords {Hamilton-Jacobi equations \and adjoint methods \and weak KAM theory \and speed of convergence.}
\subclass{35F20 \and 35F30.}
\end{abstract}
\section{Introduction}
The theory of viscosity solutions for Hamilton-Jacobi equations, introduced by Crandall and Lions in \cite{CL}, \cite{CEL} provides a body of simple and effective techniques to discover the existence, uniqueness, and stability of solutions. To date, many results concerning the speed of convergence for Hamilton-Jacobi equations of various types, which are clearly extremely important, have been studied. They are based on similar techniques,
and rely on the original ideas of Crandall and Lions. More precisely, the techniques of maximum principles, and doubling variables were used to study the speed of convergence (see \cite{BCD},  \cite{CL1}, \cite{DI}).\\\h
Recently, Evans in his forthcoming paper \cite{E} introduces some new and promising methods to study Hamilton-Jacobi equations for the time-dependent case, including the nonlinear adjoint method. This method turns out to be very useful to observe various time-dependent problems of vanishing viscosity type. The main feature of this method consists in the introduction of a further equation, to retrieve information about the solutions of the regularized problems. More precisely, one first linearizes the regularized problems, and then considers the corresponding adjoint equations. Looking at the behavior of the solutions of this new equations for different initial data, one is able to prove new estimates, which the previous techniques did not allow to obtain. In particular, the speed of convergence for time-dependent Hamilton-Jacobi equation is obtained very naturally by the adjoint method.\\\h
However, the adjoint method could not be applied directly to time-independent problems because of some difficulties such as the existence, uniqueness of the solution of the adjoint equation as well as the nonnegative property as we discuss below. In this present paper, we introduce some new ideas to apply this method to study some time-independent PDE such as stationary Hamilton-Jacobi equation, Eikonal-like equation, and Homogenization of Hamilton-Jacobi equation. We classify the problems into two classes: the class containing zero order term $u^{\epsilon}$ in the regularized equation, the class not containing zero order term $u^{\epsilon}$ in the regularized equation, and propose two different methods to deal with each class. We believe that our approach could be used in different contexts in the future.\\\h
As already mentioned, we first associate the adjoint equation to the regularized problem, and then prove an estimate for the speed of convergence. In fact, we can go further by using the constructions here and the Compensated compactness as in \cite{E1}, \cite{E} to get more properties of the solutions and get some further results. The careful study of these additional properties will be the object of a future work.\\\\\h
{\bf Outline of this paper:} This paper contains three sections about three types of static Hamilton-Jacobi equations, which are quite interesting and familiar to the readers. In order to focus on the new aspects of our approach, we keep all the hypotheses as simple as possible.\\\h
In section 2, we study the stationary problem in the whole space $\mathbb R^n$
\begin{equation}
u(x) + H(x,Du(x))=0 \h  \mbox{in} ~\mathbb R^n
\notag
\end{equation}
by looking at the regularized problem
 \begin{equation}
u^{\epsilon}(x) + H(x,Du^{\epsilon}(x))=\epsilon  \Delta u^{\epsilon} \h \mbox{in}~ \mathbb R^n.
\notag
\end{equation} 
This problem is of the first type because the regularized equation contains zero order term $u^{\epsilon}$. We have the general theme to deal with such problem like this by introducing the so-called fake parabolic adjoint equation as following:
\begin{equation}
\left\{ \begin{aligned}
-2\sigma^{\epsilon}_t - \mbox{div}(D_pH(x,Du^{\epsilon})\sigma^{\epsilon})&= \epsilon \Delta \sigma^{\epsilon} \h (x,t) \in \mathbb R^n \times (0,T)\vspace{.05in}\\
 \sigma^{\epsilon}|_T &= \delta_{x_0},\\
\end{aligned} \right. 
\notag
\end{equation}
where $\delta_{x_0}$ is the Dirac delta measure at some point $x_0 \in \mathbb R^n$.\\
The theory of distributions (see, e.g., \cite[Chapter 5]{RR}) 
ensures existence and uniqueness for $\sigma^\epsilon$, and
allows to conclude that $\sigma^\epsilon \in C^\infty(\mathbb R^n \times (0,T))$.
Studying carefully the properties of $\sigma^\epsilon$ we are able 
to obtain information about the speed of convergence of $u^{\epsilon}$
to the solution $u$. Indeed, we prove that 
$|\dfrac{\partial u^{\epsilon}}{\partial \epsilon}| \le C\epsilon^{-1/2}$,
and this in turn implies that $||u^{\epsilon} -u||_{L^{\infty}} \le C\epsilon^{1/2}$.\\\\
\h In section 3, we study the Eikonal-like equation in a bounded domain $U$ with smooth boundary
\begin{equation}
\left\{ \begin{aligned}
H(Du(x))&=0 \h  &\mbox{in}&~ U,\vspace{.05in}\\
u(x)&=0 \h  &\mbox{on}&~ \partial U.\\
\end{aligned} \right. 
\notag
\end{equation}
and we also look at the following regularized problem:
\begin{equation}
\left\{ \begin{aligned}
H(Du^{\epsilon}(x))&=\epsilon  \Delta u^{\epsilon}(x) \h  &\mbox{in}&~ U,\vspace{.05in}\\
u^{\epsilon}(x)&=0 \h  &\mbox{on}&~ \partial U.\\
\end{aligned} \right. 
\notag
\end{equation} 
This problem is of the second type since the regularized equation does not contain zero order term$u^{\epsilon}$. The idea we use to approach this type of problems is much more different to the previous one since we could not switch the problem into parabolic type. It turns out that in this case, the adjoint equation is of elliptic type and is an analog of the time-dependent Hamilton-Jacobi equation in \cite{E} as following:
\begin{equation}
\left\{ \begin{aligned}
- \mbox{div}(DH(Du^{\epsilon}) \sigma^{\epsilon})&= \epsilon  \Delta \sigma^{\epsilon}+\delta_{x_0} \h &\mbox{in}&~ U, \vspace{.05in}\\
 \sigma^{\epsilon}&=0 \h &\mbox{on}&~ \partial U.\\
\end{aligned} \right. 
\notag
\end{equation}
Besides the beauty of this adjoint equation, we furthermore can also relax the convexity condition of $H$. Up to now, all the papers dealing with the Eikonal-like equations require the convexity condition of $H$ for the bounded properties of $u^\epsilon, Du^\epsilon$ and comparison properties hence the uniqueness of the solutions. However, Ishii in \cite{I} only requires condition (H4)' instead of convexity condition for the proof of uniqueness of the solution $u$, which is the good signal for us to weaken the convexity assumption. In this section, $H$ is only required to have some kind of homogenous condition, which is much weaker, and quite natural. We will have to prove again the comparison properties and uniqueness of solutions by following the proofs in \cite{FS}, \cite{L}. Finally, we get the same speed of convergence as in the case above. One interesting point is that we could not find such result in all of the references, so it may be the new one.\\\\\h
Finally, in the last section, we study the effective Hamiltonian and homogenization of Hamilton-Jacobi equations:
\begin{equation}
H(P + Dv,y) = \bar H(P).
\notag
\end{equation}
Instead of considering the normal regularized problem as in \cite{LPV}, we consider the slightly different regularized problem, which includes the viscosity term, as following:
\begin{equation}
{\theta}z^{\theta}+ H(P+Dz^{\theta},y)= \theta^2  \Delta z^{\theta}.
\notag
\end{equation}
Following the previous classification, this is a first type problem. Thus, we introduce the fake parabolic adjoint equation as following:
\begin{equation}
\left\{ \begin{aligned}
-2\theta \sigma^{\theta}_t - \mbox{div}(D_pH\sigma^{\theta})&= \theta^2 \Delta \sigma^{\theta} \h  (x,t) \in \mathbb R^n \times (0,T)\vspace{.05in}\\
 \sigma^{\theta}|_T&=\delta_{x_0}.\\
\end{aligned} \right. 
\notag
\end{equation}
We first prove that 
$|\dfrac{\partial (\theta z^{\theta})}{{\partial \theta}}| \le C$ and then, 
as a consequence, we get the estimate
$||\theta z^{\theta}+ \bar H(P)||_{L^{\infty}} \le C\theta$. 
This result, in particular, implies what obtained in \cite{DI}.
\section{Stationary problem in $\mathbb R^n$}
We are going to study the stationary Hamilton-Jacobi equation in $\mathbb R^n$
\begin{equation}
u(x) + H(x,Du(x))=0 \h \mbox{in}~ \mathbb R^n.
\label{HJ1}
\end{equation} 
As usual, we consider the following regularized equation
\begin{equation}
u^{\epsilon}(x) + H(x,Du^{\epsilon}(x))=\epsilon  \Delta u^{\epsilon} \h \mbox{in}~ \mathbb R^n.
\label{HJ2}
\end{equation} 
Let us for simplicity assume that $H$ is smooth and $H$ satisfies some conditions as in \cite{L}, \cite{S} as following
\begin{equation}
\left\{ \begin{aligned}
\sup_{x \in \mathbb R^n}|H(x,0)| \le C \le \infty; ~\sup_{x \in \mathbb R^n}|D_xH(x,p)| \le C(1+|p|),\vspace{.05in}\\
H(x,p) \to \infty ~\mbox{as}~ |p| \to \infty ~ \mbox{uniformly in} ~x \in \mathbb R^n.\hh\\
\end{aligned} \right. 
\notag
\end{equation}
By the coercive property of the Hamiltonian $H$ and the properties above, we have some well-known standard estimates from \cite{L} as following:
\begin{equation}
 ||u^{\epsilon}||_{L^{\infty}}, ||Du^{\epsilon}||_{L^{\infty}} \le C.
 \label{HJ3}
 \end{equation}
Our main theorem of this section is
\begin{Theorem}
There exists a constant $C > 0$ such that
\begin{equation}
||u^{\epsilon}-u||_{L^{\infty} }\le C \epsilon^{1/2}.
\label{HJ4}
\end{equation} 
\label{THJ1}
\end{Theorem}
In fact, this theorem was proved long time ago, for instance in \cite{CL1}, \cite{F}, \cite{S}, \cite{BCD}. However, we propose here a new way to prove it by using adjoint method.
\begin{Lemma}
Let $w^{\epsilon} = \dfrac{|Du^{\epsilon}|^2}{2}$ then $w^{\epsilon}$ satisfies:
\begin{equation}
2w^{\epsilon}+D_{p}H(x,Du^{\epsilon}).Dw^{\epsilon} +D_{x}H(x,Du^{\epsilon}).Du^{\epsilon}= \epsilon\Delta w^{\epsilon} - \epsilon |D^2u^{\epsilon}|^2.
\label{HJ5}
\end{equation}
\label{LHJ1}
\end{Lemma}
{\bf Proof}\\
Differentiate the equation (\ref{HJ2}) with respect to $x_i$
\begin{equation}
 u^{\epsilon}_{x_i} + H_{x_i}(x, Du^{\epsilon})+H_{p_k}(x,Du^{\epsilon})u^{\epsilon}_{x_k x_i} = \epsilon \Delta u^{\epsilon}_{x_i}.
\label{HJ6}
\end{equation}
Taking the product of (\ref{HJ6}) with $u^\epsilon_{x_i}$ and summing over $i$
\begin{equation}
 |Du^{\epsilon}|^2 + D_{x}H(x, Du^{\epsilon}).Du^{\epsilon}+H_{p_k}(x,Du^{\epsilon})(\dfrac{|Du^{\epsilon}|^2}{2})_{x_k} = \epsilon \Delta u^{\epsilon}_{x_i}u^{\epsilon}_{x_i}.
\label{HJ7}
\end{equation}
Furthermore, notice that:
\begin{equation}
\Delta u^{\epsilon}_{x_i}u^{\epsilon}_{x_i} = u^{\epsilon}_{x_k x_k x_i}u^{\epsilon}_{x_i} =  (u^{\epsilon}_{x_k x_i}u^{\epsilon}_{x_i})_{x_k}  - \sum_{i,k} | u^{\epsilon}_{x_k x_i}|^2 = \Delta (\dfrac{|Du^{\epsilon}|^2}{2}) - |D^2u^{\epsilon}|^2.
\notag
\end{equation}
Combining those two calculations, we get the lemma.\\\\
Now we introduce the new function $v^{\epsilon}$ to change ($\ref{HJ5}$) into a "fake parabolic" equation.  We will explain later the reason why we have to switch to parabolic type. Let $T>0$ be a constant and let  \vspace{.05in}\\\hhh\hh
$v^{\epsilon} (x,t) = e^t w^{\epsilon}(x) \h (x,t) \in \mathbb R^n \times [0,T]$. \vspace{.05in}
\\ Then from ($\ref{HJ5}$), we therefore get that $v^{\epsilon}$ satisfies:
\begin{equation}
2v^{\epsilon}_t+D_{p}H(x,Du^{\epsilon}).Dv^{\epsilon} +e^tD_{x}H(x,Du^{\epsilon}).Du^{\epsilon}= \epsilon\Delta v^{\epsilon} - \epsilon e^t|D^2u^{\epsilon}|^2.
\label{HJ8}
\end{equation}
\\
{\bf Adjoint method:} We now introduce the adjoint equation of equation (\ref{HJ8}). For $x_0 \in \mathbb R^n$, let $\sigma^{\epsilon}$ be the solution of the following PDE:
\begin{equation}
\left\{ \begin{aligned}
-2\sigma^{\epsilon}_t - \mbox{div}(D_pH(x,Du^{\epsilon})\sigma^{\epsilon})&= \epsilon \Delta \sigma^{\epsilon} \h  (x,t) \in \mathbb R^n \times (0,T)\vspace{.05in}\\
 \sigma^{\epsilon}|_T&=\delta_{x_0}.\\
\end{aligned} \right. 
\label{HJ9}
\end{equation}
By the theory of distributions (see Chapter 5 in \cite{RR}), $\sigma^\epsilon$ is unique and $\sigma^\epsilon \in C^\infty(\mathbb R^n \times (0,T))$.\\
From the solution $\sigma^{\epsilon}$ of the adjoint equation, we can somehow figure out the properties of $u^{\epsilon}$ as well as $u$, which are our very important goals especially in the case that $H$ is not convex in $p$.\\ 
Firstly, let us point out some properties of $\sigma^{\epsilon}$:
\begin{Lemma} {\bf Properties of $\sigma^{\epsilon}$}
\begin{itemize}
\item[(i)] $\sigma^{\epsilon}(x,t) \ge 0$~ for $(x,t) \in \mathbb R^n \times (0,T)$,\\
 \item[(ii)] $\int_{\mathbb R^n} \sigma^{\epsilon}(x,t)dx=1$~ for $ t \in (0,T)$.
 \end{itemize}
\label{LHJ2}
\end{Lemma}
{\bf Proof}\\
The proof can be easily obtained by using the Maximum Principle and by integrating over $\mathbb R^n$.\\ 
\begin{Remark}
\label{Adjoint1}
\end{Remark}
As we can see, when we switch the equation into the parabolic type then we automatically have the existence of the solution $\sigma^{\epsilon}$ of the adjoint equation as well as the maximum principle can be applied with the only requirement of the boundedness of coefficients. Note that we need the property (i) of the above Lemma to do further derivations as you can see below.\\
Besides, one can write down the adjoint equation of (\ref{HJ5}) in form of elliptic equation and can see that the adjoint equation may not have the solution, the uniqueness of the solution as well as the required condition to apply the maximum principle.\\\\
Now, we start to observe properties and connections between $\sigma^{\epsilon}$ and $u^{\epsilon}$
\begin{Lemma}
There exists a constant $C>0$ such that
\begin{equation}
\int_{0}^{T} \int_{\mathbb R^n} \epsilon e^t|D^2u^{\epsilon}|^2 \sigma^{\epsilon}dxdt \le C.
\label{HJ10}
\end{equation}
\label{LHJ3}
\end{Lemma}
{\bf Proof}\\
\begin{align}
\label{HJ11}
&\dfrac{d}{dt}\int_{ \mathbb R^n} 2 \sigma^{\epsilon} v^{\epsilon}=\int_{\mathbb R^n}2 \sigma^{\epsilon}_t v^{\epsilon}+2 \sigma^{\epsilon} v^{\epsilon}_t \qquad \qquad \qquad \qquad \qquad \qquad \qquad \qquad \qquad
\\
=&\int_{\mathbb R^n}2 \sigma^{\epsilon}_t v^{\epsilon}+\int_{\mathbb R^n}(-D_{p}H.Dv^{\epsilon} -e^t D_{x}H.Du^{\epsilon}+ \epsilon \Delta v^{\epsilon} - \epsilon  e^t  |D^2u^{\epsilon}|^2) \sigma^{\epsilon} \quad\qquad
\notag\\
=&\int_{\mathbb R^n}(2 \sigma^{\epsilon}_t+\mbox{div}(D_pH \sigma^{\epsilon})+\epsilon \Delta \sigma^{\epsilon}) v^{\epsilon}-\int_{\mathbb R^n}(e^t D_{x}H.Du^{\epsilon}+\epsilon  e^t  |D^2u^{\epsilon}|^2) \sigma^{\epsilon}\qquad
\notag\\
=&-\int_{\mathbb R^n}(e^t D_{x}H.Du^{\epsilon}+\epsilon e^t  |D^2u^{\epsilon}|^2) \sigma^{\epsilon}. \qquad \qquad \qquad \qquad \qquad \qquad \qquad \qquad\quad
\notag
\end{align}
Now we integrate (\ref{HJ11}) from $0$ to $T$:
\begin{align}
\label{HJ12}
&\int_{\mathbb R^n} 2 \sigma^{\epsilon}(x,T) v^{\epsilon}(x,T)dx-\int_{\mathbb R^n} 2 \sigma^{\epsilon}(x,0) v^{\epsilon}(x,0)dx \qquad \qquad \qquad \qquad \qquad 
\\
=&-\int_{0}^{T}\int_{\mathbb R^n}e^t D_{x}H.Du^{\epsilon} \sigma^{\epsilon} -\int_{0}^{T}\int_{\mathbb R^n}\epsilon  e^t  |D^2u^{\epsilon}|^2 \sigma^{\epsilon}dxdt.\qquad \qquad \qquad \qquad 
\notag
\end{align}
Hence we get:
\begin{align}
\label{HJ13}
&\int_{0}^{T}\int_{\mathbb R^n}\epsilon e^t |D^2u^{\epsilon}|^2 \sigma^{\epsilon}dxdt \qquad \qquad \qquad \qquad \qquad \qquad \qquad \qquad \qquad\\
 \le& ~| 2 v^{\epsilon}(x_0,T)|+
|\int_{\mathbb R^n} 2 \sigma^{\epsilon}(x,0) v^{\epsilon}(x,0)dx|+|\int_{0}^{T}\int_{\mathbb R^n}e^t D_{x}H.Du^{\epsilon} \sigma^{\epsilon}| \quad
\notag\\
\le& ~2 e^T C + 2 C + C(e^T-1) \le C. \qquad \qquad \qquad \qquad \qquad \qquad \qquad \qquad
\notag
\end{align}
We get the lemma. As stated in \cite{E}, the estimate (\ref{HJ10}) seems to be the new estimate and it will be very helpful later.\\\\
Notice that all of the estimates here are independent of the choice of $x_0$. More precisely, for any $x_0 \in \mathbb R^n$ and the corresponding $\sigma^{\epsilon}$, the estimates stay the same with the same constants. More generally, we also have all such estimates if we assume $\sigma^{\epsilon}|_T = \nu$ for $\nu$ is a Borel probability measure. However, we do not really use the general probability measure here.
\begin{Definition}
Define $u_{\epsilon}^{\epsilon}(x) = \dfrac{\partial u^{\epsilon}}{\partial \epsilon}(x)$.
\label{DHJ1}
\end{Definition}
We have the following theorem:
\begin{Theorem}
There exists a constant $C>0$ such that
\begin{equation}
|u_{\epsilon}^{\epsilon}(x)| \le C  {\epsilon}^{-1/2}.
\label{HJ14}
\end{equation}
\label{THJ2}
\end{Theorem}
{\bf Proof} \\
According to standard elliptic estimates, the function $u^{\epsilon}$ is smooth in the parameter $\epsilon$ away from $\epsilon=0$. Differentiate (\ref{HJ2}) with respect to ${\epsilon}$
\begin{equation}
u_{\epsilon}^{\epsilon} + D_pH(x,Du^{\epsilon}).Du_{\epsilon}^{\epsilon} = \epsilon \Delta u_{\epsilon}^{\epsilon} + \Delta u^{\epsilon}.
\label{HJ15}
\end{equation}
Define $z^{\epsilon}:\mathbb R^n \times [0,T] \to \mathbb R$ such that $z^{\epsilon}(x,t)=e^t u_{\epsilon}^{\epsilon}(x)$. Then $z^{\epsilon}$ satisfies the following PDE:
\begin{equation}
z^{\epsilon}_t + D_pH(x,Du^{\epsilon}).Dz^{\epsilon} = \epsilon  \Delta z^{\epsilon} + e^t \Delta u^{\epsilon}.
\label{HJ16}
\end{equation}
Notice that the coefficients of (\ref{HJ16}) is slightly different to those of the adjoint equation. Playing the same tricks as in Lemma \ref{LHJ3} we have:
\begin{align}
\label{HJ17}
&\dfrac{d}{dt}\int_{\mathbb R^n} 2 \sigma^{\epsilon} z^{\epsilon}=\int_{\mathbb R^n}2 \sigma^{\epsilon}_t v^{\epsilon}+\sigma^{\epsilon} z^{\epsilon}_t +\int_{\mathbb R^n} e^t u_{\epsilon}^{\epsilon} \sigma^{\epsilon}\qquad \qquad \qquad \qquad \qquad \qquad
\\
=&\int_{\mathbb R^n}2 \sigma^{\epsilon}_t z^{\epsilon}+\int_{\mathbb R^n}(-D_{p}H.Dz^{\epsilon}+ \epsilon \Delta z^{\epsilon}+e^t \Delta u^{\epsilon}) \sigma^{\epsilon} +\int_{\mathbb R^n} e^t u_{\epsilon}^{\epsilon} \sigma^{\epsilon} \qquad \quad\qquad
\notag\\
=&\int_{\mathbb R^n}(2 \sigma^{\epsilon}_t+\mbox{div}(D_pH \sigma^{\epsilon})+\epsilon \Delta \sigma^{\epsilon}) z^{\epsilon}+\int_{\mathbb R^n}(e^t \Delta u^{\epsilon} \sigma^{\epsilon}+e^t u_{\epsilon}^{\epsilon} \sigma^{\epsilon}) \qquad \qquad \qquad
\notag\\
=&\int_{\mathbb R^n}(e^t \Delta u^{\epsilon} \sigma^{\epsilon}+e^t u_{\epsilon}^{\epsilon} \sigma^{\epsilon}). \qquad \qquad \qquad \qquad \qquad \qquad \qquad \qquad \qquad \qquad \qquad
\notag
\end{align}
By appropriate pertubation arguments which are quite usual and classical in the use of maximum principle, for example in \cite{LM}, we can assume without loss of generality that there exists $x_1 \in \mathbb R^n$ such that 
$$
|u_{\epsilon}^{\epsilon}(x_1)| = \max_{\mathbb R^n} |u_{\epsilon}^{\epsilon}(x)|.
$$
Now we let $x_0=x_1$ in the adjoint equation (\ref{HJ9}).\\
Integrate (\ref{HJ17}) from $0$ to $T$
\begin{align}
\label{HJ18}
&\int_{\mathbb R^n} 2 \sigma^{\epsilon}(x,T) z^{\epsilon}(x,T)dx-\int_{\mathbb R^n} 2 \sigma^{\epsilon}(x,0) z^{\epsilon}(x,0)dx \qquad \qquad \qquad \qquad \qquad 
\\
=&\int_{0}^{T}\int_{\mathbb R^n}e^t \Delta u^{\epsilon} \sigma^{\epsilon}+\int_{0}^{T}\int_{\mathbb R^n}e^t u_{\epsilon}^{\epsilon} \sigma^{\epsilon}dxdt.\qquad \qquad \qquad \qquad \qquad \qquad \qquad \quad
\notag
\end{align}
Substitute the condition $\sigma^{\epsilon}|_T = \delta_{x_0}$ into the equation above,
\begin{equation}
|2 e^T u_{\epsilon}^{\epsilon}(x_0) - \int_{\mathbb R^n} 2 \sigma^{\epsilon}(x,0) z^{\epsilon}(x,0) - \int_{0}^{T}\int_{\mathbb R^n}e^t u_{\epsilon}^{\epsilon} \sigma^{\epsilon}|=|\int_{0}^{T}\int_{\mathbb R^n}e^t \Delta u^{\epsilon} \sigma^{\epsilon}|.
\label{HJ19}
\end{equation}
Furthermore, we can control the left hand side of (\ref{HJ19}) as following
\begin{align}
\mbox{LHS}&=|2 e^T u_{\epsilon}^{\epsilon}(x_0) - \int_{\mathbb R^n} 2 \sigma^{\epsilon}(x,0) z^{\epsilon}(x,0) - \int_{0}^{T}\int_{\mathbb R^n}e^t u_{\epsilon}^{\epsilon} \sigma^{\epsilon}| \qquad \qquad \quad\notag\\
&\ge 2 e^T |u_{\epsilon}^{\epsilon}(x_0)|- \int_{\mathbb R^n} 2 \sigma^{\epsilon}(x,0) |u_{\epsilon}^{\epsilon}(x_0)| - \int_{0}^{T}\int_{\mathbb R^n}e^t |u_{\epsilon}^{\epsilon}(x_0)| \sigma^{\epsilon} \qquad \quad\notag\\
&=|u_{\epsilon}^{\epsilon}(x_0)|  (2 e^T -2 - (e^T-1)) = |u_{\epsilon}^{\epsilon}(x_0)|  (e^T-1).\qquad \qquad \qquad \quad \notag
\end{align}
Besides, by using Lemma \ref{LHJ3} and Holder's inequality, we can estimate the right hand side of (\ref{HJ19}):
\begin{align}
\mbox{RHS} &= |\int_{0}^{T}\int_{\mathbb R^n}e^t \Delta u^{\epsilon} \sigma^{\epsilon}| \le C \int_{0}^{T}\int_{\mathbb R^n}|D^2 u^{\epsilon}| \sigma^{\epsilon}\qquad \qquad \qquad\notag \\
&\le C  \{ \int_{0}^{T}\int_{\mathbb R^n}|D^2 u^{\epsilon}|^2 \sigma^{\epsilon}\}^{1/2}  \{ \int_{0}^{T}\int_{\mathbb R^n}\sigma^{\epsilon}\}^{1/2} \le C  {\epsilon}^{-1/2}.\qquad\notag
\end{align}
So we get the theorem.\\\\
{\bf Proof of Theorem \ref{THJ1}}\\
By using Theorem \ref{THJ2}, we immediately get the result.
\section{Eikonal-like equation in bounded domain}
We are going to study the following of Eikonal-like Hamilton-Jacobi equation in the given bounded domain $U$ with smooth boundary:
\begin{equation}
\left\{ \begin{aligned}
H(Du(x))&=0 \h  &\mbox{in}&~ U,\vspace{.05in}\\
u(x)&=0 \h &\mbox{on}&~ \partial U.\\
\end{aligned} \right. 
\label{HJ20}
\end{equation}
Our approach, as usual, is to consider regularized problem:
\begin{equation}
\left\{ \begin{aligned}
H(Du^{\epsilon}(x))&=\epsilon   \Delta u^{\epsilon}(x) \h &\mbox{in}&~ U,\vspace{.05in}\\
u^{\epsilon}(x)&=0 \h &\mbox{on}&~ \partial U.\\
\end{aligned} \right. 
\label{HJ21}
\end{equation} 
\h Crandall and Lions study this equation in sense of viscosity solution first in \cite{CL} and Lions also studies it in \cite{L}. After that, Fleming and Souganidis study it in more details and also give some asymptotic series of the solutions of the regularized problem in \cite{FS}. Then Ishii gives a simple and direct proof of the uniqueness of the solution in \cite{I}. We here base on the conditions given in \cite{FS}, \cite{L} and we refer the readers to \cite{I}, \cite{FS} and \cite{L} for more details.\\\h
Our goal here is not only to prove the speed of convergence but also to relax the convexity conditions of $H$. Obviously we cannot relax the convexity condition without require some sufficient conditions as we will see in the counter-example below. But the condition we need is much more weaker and quite natural like the homogenous condition. We assume $H$ satisfies the following conditions:\\\\\h
(H1) $H$ smooth and $H(0) < 0$,\\\h
(H2) $H$ is superlinear, i.e. $\lim_{|p| \to \infty} \dfrac{H(p)}{|p|}=\infty$,\\\h
(H3) There exist $\gamma, \delta > 0$ s.t. $DH(p).p - \gamma H(p) \ge \delta >0 \h \forall~ p \in \mathbb R^n$.\\\\
The condition (H3) is used to replace the convexity condition and will be discussed later. We just make an obvious observation that if $H$ is convex then we have (H3) with $\gamma=1$ and $\delta=-H(0)$.\\
\begin{Theorem}
There exists $C>0$ such that
\begin{equation}
||u^{\epsilon}||_{L^{\infty}}, ||Du^{\epsilon}||_{L^{\infty}} \le C.  
 \label{HJ22}
 \end{equation}
 \label{FS1}
 \end{Theorem}
 {\bf Proof}\\
In the case where $H$ is convex then this theorem is proved in \cite{FS} by Lemma 1.1 and 1.2 or in \cite{L}. Here, we follow almost all of the proofs and just need to slightly change some estimates that use the convexity condition.\\
By Lemma 1.1 and the first part of Lemma 1.2 in \cite{FS}, there exists a constant $C>0$ such that
$0 \le u^{\epsilon} \le C$ in $\bar U$ and $|Du^{\epsilon}| \le C$ on $\partial U$.\\
To complete the proof, we will only need to bound $|Du^{\epsilon}|$ in $U$.\\
Using the same ideas like  in \cite{FS}, \cite{L}, let $w = |Du^{\epsilon}| - \mu u^{\epsilon}$, where $\mu$ is to be a suitably chosen constant. Suppose that $w$ has a positive maximum at an interior point $x_0 \in U$. At $x_0$ we have:
 \begin{equation}
0=w_{x_i} = \dfrac{\sum_{k}u_{x_k}^{\epsilon}u_{x_k x_i}^{\epsilon}}{|Du^{\epsilon}|} - \mu u_{x_i}^{\epsilon},  
 \notag
 \end{equation}
Hence we get:
\begin{equation}
\sum_{i} (\sum_{k}u_{x_k}^{\epsilon}u_{x_k x_i}^{\epsilon})^2 = \mu^2 |Du^{\epsilon}|^4.
\notag
\end{equation}
Furthermore, we also have:
 \begin{equation}
0  \le -\epsilon \Delta w =  \dfrac{\epsilon \sum_{i} (\sum_{k}u_{x_k}^{\epsilon}u_{x_k x_i}^{\epsilon})^2}{|Du^{\epsilon}|^3} - \dfrac{\epsilon \sum_{i,k} (u_{x_k x_i}^{\epsilon})^2}{|Du^{\epsilon}|}+\dfrac{ \sum_{k} u_{x_k}^{\epsilon}(-\epsilon \Delta u^{\epsilon})_{x_k}}{|Du^{\epsilon}|}+\mu (\epsilon \Delta u^{\epsilon}),
\notag
\end{equation}
By using the inequality $\dfrac{(\Delta u^{\epsilon})^2}{n} \le \sum_{i,k} (u^{\epsilon}_{x_i x_k})^2$ and (\ref{HJ21})
\begin{equation}
0 \le \epsilon \mu^2 |Du^{\epsilon}| - \dfrac{H^2}{n \epsilon |Du^{\epsilon}|} - \mu DH.Du^{\epsilon} + \mu H.
\notag
\end{equation}
 Besides, (H3) implies
 \begin{equation}
\mu DH.Du^{\epsilon} - \mu \gamma H > \delta \mu>0,
\notag
\end{equation}
Thus,
\begin{equation}
\dfrac{H^2}{ |Du^{\epsilon}|^2} \le n \mu^2 \epsilon^2 + n \epsilon (\mu - \mu \gamma)\dfrac{H}{|Du^{\epsilon}|}  \le n \mu^2 \epsilon^2 + n \epsilon \mu (1+\gamma) \dfrac{|H|}{|Du^{\epsilon}|}.
\notag
\end{equation}
Choose $\mu=\dfrac{1}{2n(1+\gamma)}$ then for $\epsilon <1$, we get the estimate:
\begin{equation}
\dfrac{H^2}{ |Du^{\epsilon}|^2} \le 1+ \dfrac{|H|}{|Du^{\epsilon}|}.
\notag
\end{equation}
By the superlinearity condition (H2) we finally get $|Du^{\epsilon}|$ is bounded independently of $\epsilon$. We get the theorem.
\begin{Remark}
\label{RHJ1}
\end{Remark}
The existence of the solution of (\ref{HJ21}) then follows directly from \cite{FS} with some changes and adaptations similar to the proof of Theorem \ref{FS1} above.\\
Now we discuss about the uniqueness of the viscosity solution $u$ of (\ref{HJ20}).\\
For $p \in \mathbb R^n$, let's consider the following function $\phi$ from $(0,\infty)$ to $\mathbb R$
\begin{equation}
\phi(t) = t^{-\gamma}H(tp) \quad \forall~ t>0,
\notag
\end{equation}
then
\begin{equation}
\phi'(t) = t^{-\gamma-1}(DH(tp).(tp) - \gamma H(tp)) >  t^{-\gamma-1}\delta >0.
\notag
\end{equation}
Hence $\phi$ is strictly increasing and for $t<1$ we have furthermore:
\begin{equation}
\phi(1) - \phi(t) = \int_{t}^1 \phi'(s)ds >\int_{t}^1 s^{-\gamma-1}\delta ds =\dfrac{\delta}{\gamma+1} (t^{-\gamma} -1) >0,
\notag
\end{equation}
Thus,
\begin{equation}
H(tp) \le t^{\gamma}H(p) - \dfrac{\delta}{\gamma+1} (1-t^{\gamma})=t^{\gamma}H(p) + \dfrac{-\delta}{(\gamma+1)H(0)} (1-t^{\gamma})H(0).
\notag
\end{equation}
By (H1) we have that $H(0)<0$. So $H$ satisfies all the conditions (H1)-(H3) and (H4)' of \cite{I} with $\varphi=0$. Hence the uniqueness of viscosity solution of (\ref{HJ20}) follows.\\
The proof of the uniqueness of $u^{\epsilon}$ is quite complicated and follows the key idea of this section. Therefore, we put it in the appendix at the end of this paper.\\\\
Our main theorem of this section is
\begin{Theorem}
There exists a constant $C > 0$ such that
\begin{equation}
||u^{\epsilon}-u||_{L^{\infty} }\le C  \epsilon^{1/2}.
\label{HJ23}
\end{equation} 
\label{THJ3}
\end{Theorem}
Some of the lemmas below will be quite similar to those in Section 2. Therefore, they will only be stated without proofs unless there are some huge differences.
\begin{Lemma}
Let $w^{\epsilon} = \dfrac{|Du^{\epsilon}|^2}{2}$ then $w^{\epsilon}$ satisfies:
\begin{equation}
DH(Du^{\epsilon}).Dw^{\epsilon} = \epsilon \Delta w^{\epsilon} - \epsilon   |D^2u^{\epsilon}|^2.
\label{HJ24}
\end{equation}
\label{LHJ4}
\end{Lemma}
Note that the term $w^{\epsilon}$ does not appear in equation (\ref{HJ24}). Hence we cannot convert this equation to the parabolic type as in section 2. We then introduce a different approach.\\\\
{\bf Adjoint method:} We now introduce the adjoint equation to the equation (\ref{HJ24}). For each $x_0 \in U$, we consider the following PDE:
\begin{equation}
\left\{ \begin{aligned}
- \mbox{div}(DH(Du^{\epsilon}) \sigma^{\epsilon})&= \epsilon  \Delta \sigma^{\epsilon}+\delta_{x_0} \h &\mbox{in}&~ U, \vspace{.05in}\\
 \sigma^{\epsilon}&=0 \h &\mbox{on}&~ \partial U.\\
\end{aligned} \right. 
\label{HJ25}
\end{equation}
The adjoint equation here is very nice and somehow similar to the one that Evans introduces in \cite{E} for the time-dependent case. From $\sigma^{\epsilon}$, the solution of the adjoint equation, we can somehow figure out the properties of $u^{\epsilon}$ as well as $u$, which are our very important goals especially in the case that $H$ is not convex in $p$. However, the problem is, like what we have mentioned in the Remark \ref{Adjoint1} above, we do not know about the existence, uniqueness of (\ref{HJ25}) as well as the nonnegative property of $\sigma^{\epsilon}$, which we really need.\\ 
It's quite interesting that to observe $\sigma^{\epsilon}$, we once again need the adjoint equation of (\ref{HJ25}):\\\h
For each $f \in C^\infty(\bar U)$ and $ f \ge 0$, we consider the following equation
\begin{equation}
\left\{ \begin{aligned}
DH(Du^{\epsilon}).Dv^{\epsilon}&= \epsilon \Delta v^{\epsilon}+f \h &\mbox{in}&~ U, \vspace{.05in}\\
v^{\epsilon}&=0 \h &\mbox{on}&~ \partial U.\\
\end{aligned} \right. 
\label{HJ26}
\end{equation}
For $f=0$ then it's obvious by the Maximum principle that $v^{\epsilon}=0$. Hence by Fredholm alternative, both equations (\ref{HJ26}) and (\ref{HJ25}) have unique solutions. By the theory of distributions (see Chapter 5 in \cite{RR}), $\sigma^\epsilon \in C^\infty(U\setminus\{x_0\})$.\\
Furthermore, by Maximum principle again, $v^{\epsilon} \ge 0$.
\begin{Lemma}
The following fact holds  
\begin{equation}
\int_{U} f \sigma^{\epsilon} dx = v^{\epsilon}(x_0) \ge 0.
\label{HJ27}
\end{equation}
Hence in particular, $\sigma^{\epsilon} \ge 0$ in $U \setminus \{x_0\}$.
\label{LHJ5}
\end{Lemma}
{\bf Proof}\\
By (\ref{HJ25}) and (\ref{HJ26})
\begin{align}
\label{HJ28}
&\int_{U} f \sigma^{\epsilon} dx=\int_{U} DH(Du^{\epsilon}).Dv^{\epsilon} \sigma^{\epsilon}-{\epsilon} \Delta v^{\epsilon} \sigma^{\epsilon}\qquad \qquad \qquad \qquad \qquad 
\\
=&\int_{U} (-\mbox{div}(DH(Du^{\epsilon}) \sigma^{\epsilon})-\epsilon \Delta \sigma^{\epsilon}) v^{\epsilon}=v^{\epsilon}(x_0) \ge 0.\qquad \qquad \qquad \qquad
\notag
\end{align}
We therefore get the lemma.\\
From the above lemma, we can easily derive some following properties of $\sigma^{\epsilon}$:
\begin{Lemma} {\bf Properties of $\sigma^{\epsilon}$}
\begin{itemize}
\item[(i)] $ \sigma^{\epsilon}\ge 0$ in $U \setminus \{x_0\}$. In particular, $\dfrac{\partial \sigma^{\epsilon}}{\partial n} \le 0$ on $\partial U$.\\
 \item[(ii)] $\int_{\partial U}\epsilon  \dfrac{\partial \sigma^{\epsilon}}{\partial n} dS=-1$.
 \end{itemize}
\label{LHJ6}
\end{Lemma}

\begin{Lemma}
There exists a constant $C>0$ such that
\begin{equation}
\int_{U} \epsilon |D^2u^{\epsilon}|^2 \sigma^{\epsilon}dx\le C.
\label{HJ29}
\end{equation}
\label{LHJ7}
\end{Lemma}
{\bf Proof}\\
By (\ref{HJ24}), we have:
\begin{equation}
\int_{U} (DH(Du^{\epsilon}).Dw^{\epsilon}-\epsilon \Delta w^{\epsilon}) \sigma^{\epsilon}dx=-\int_{U} \epsilon |D^2u^{\epsilon}|^2 \sigma^{\epsilon}dx.
\label{HJ30}
\end{equation}
Integrate by parts the left hand side of the above equality:
\begin{align}
\label{HJ31}
\mbox{LHS}&=\int_{U} -\mbox{div}(DH(Du^{\epsilon}) \sigma^{\epsilon}) w^{\epsilon}-\epsilon \Delta \sigma^{\epsilon} w^{\epsilon}+\int_{\partial U} \epsilon  \dfrac{\partial \sigma^{\epsilon}}{\partial n} w^{\epsilon} \qquad \qquad \qquad 
\\
&=\int_{U} (-\mbox{div}(DH(Du^{\epsilon}) \sigma^{\epsilon})-\epsilon \Delta \sigma^{\epsilon}) w^{\epsilon}+\int_{\partial U} \epsilon  \dfrac{\partial \sigma^{\epsilon}}{\partial n} w^{\epsilon}\qquad \qquad \qquad \quad
\notag \\
&=w(x_0)  +\int_{\partial U} \epsilon  \dfrac{\partial \sigma^{\epsilon}}{\partial n} w^{\epsilon}.\qquad \qquad \qquad \qquad \qquad\qquad\qquad\qquad\qquad\quad\notag
\end{align}
So, by using Lemma \ref{LHJ6}, we get the lemma.\\\\
As normal, if we can bound $\int_{U} \sigma^{\epsilon}dx$ independently of $\epsilon$ then the result follows immediately as one can see later by using the same arguments as in Section 2. However, it's not easy to bound $\int_{U} \sigma^{\epsilon}dx$ here. We will show the reasons why in the following discussions.\\
Choose $f=1$ then ($\ref{HJ26}$) reads
\begin{equation}
\left\{ \begin{aligned}
DH(Du^{\epsilon}).Dv^{\epsilon}&= \epsilon  \Delta v^{\epsilon}+1 \h &\mbox{in}&~ U, \vspace{.05in}\\
v^{\epsilon}&=0 \h &\mbox{on}&~ \partial U.\\
\end{aligned} \right. 
\label{HJ32}
\end{equation}
And also Lemma $\ref{LHJ5}$ reads
\begin{equation}
\int_{U} \sigma^{\epsilon} dx = v^{\epsilon}(x_0) \ge 0.
\label{HJ33}
\end{equation}
Hence, in order to bound $\int_{U} \sigma^{\epsilon}dx$, we need to bound $ v^{\epsilon}(x_0)$. And since $x_0$ may vary, $\max_{U} v^{\epsilon}$ should be bounded uniformly independently of $\epsilon$.\\
It turns out that this fact is not true for general $H$. For example, when $DH(p)=0$ for all $p$, the above fact is no longer true, i.e. we will no longer have the uniformly bound for $\max_{U} v^{\epsilon}$ by the following explicit example.\\
Let's consider the following ODE:
\begin{equation}
\left\{ \begin{aligned}
\epsilon \Delta v^{\epsilon}+1&=0 \hh \mbox{in}~ (0,1), \vspace{.05in}\\
v^{\epsilon}(0)=v^{\epsilon}(1)&=0. \\
\end{aligned} \right. 
\label{HJ34}
\end{equation}
Then $v^{\epsilon}(x)= \dfrac{1}{2 \epsilon}(x-x^2)$, which implies $\max_{[0,1]} v^{\epsilon}= \dfrac{1}{8 \epsilon}$. So $\max_{[0,1]} v^{\epsilon}$ blows up as $\epsilon$ tends to $0$.\\
Heuristically, this counter-example shows that we need to have some conditions on the gradient of the Hamiltonian $H$ that allow us to control $v^\epsilon$.\\
We introduce next the second example, where we have some growth control on $DH(p)$, as following
\begin{equation}
\left\{ \begin{aligned}
(v^{\epsilon})'&= \epsilon \Delta v^{\epsilon}+1 \hh \mbox{in}~ (0,1), \vspace{.05in}\\
v^{\epsilon}(0)&=v^{\epsilon}(1)=0. \\
\end{aligned} \right. 
\label{HJ35}
\end{equation}
Then we get
\begin{equation}
v^{\epsilon} (x )= x-\dfrac{e^{x/{\epsilon}}-1}{e^{1/{\epsilon}}-1}.
\notag
\end{equation}
Then $\max_{[0,1]} v^{\epsilon} \le 1$, which provides us the uniformly boundedness of $\max_{U} v^{\epsilon}$ independent of $\epsilon$.\\
While the first example fails, the second one intuitively shows that with good growth control on $DH(p)$, we will have such uniform bound.\\\\
Based upon the above examples and discussions, we introduce the following condition (H3) to have the uniform bound of $\max_U v^\epsilon$ independent of $\epsilon$
\vspace{.05in}\\\h
(H3) There exist $\gamma, \delta > 0$ s.t. $DH(p).p - \gamma H(p) \ge \delta >0 \h \forall~ p \in \mathbb R^n$. \vspace{.05in}\\
In particular, if we choose $\gamma=1, \delta=-H(0)$, then (H3) becomes
\begin{equation}
DH(p).p - H(p) \ge -H(0) >0 \h \forall ~p \in \mathbb R^n,
\notag
\end{equation}
which is the convexity-like condition for $H$. And also if $H$ is convex then (H3) follows with $\gamma=1, \delta=-H(0)$.\\
In fact, the required condition (H3) is similar to the homogenous condition. It's natural and it works well for a lot of cases where $H$ is not convex. For example, for $n=1$, let's consider the following function:
\begin{equation}
H(p)=(p^2-1)^2-2=p^4-2p^2-1,
\notag
\end{equation}
then $H$ is not convex and
\begin{equation}
DH(p).p-2H(p)=(4p^4-4p^2)-2(p^4-2p^2-1)=2p^4+2 \ge 2>0.
\notag
\end{equation}
It's easy to see that $H$ satisfies all the required conditions of our problem even though $H$ is not convex.\\
Furthermore, the condition (H3) is suitable and fit well in every required step of our problem.\\ 
The following lemma is the key lemma of this section, it shows the way to bound $\max_{U} v^{\epsilon}$: 
\begin{Lemma}
Let $\alpha, \beta \in \mathbb R$ and $z(x) = \alpha x.Du^{\epsilon}(x) + \beta u^{\epsilon}(x)$ then
\begin{equation}
DH(Du^{\epsilon}).Dz - \epsilon  \Delta z = (\alpha+\beta) DH(Du^{\epsilon}).Du^{\epsilon} - (2 \alpha+\beta) \epsilon \Delta u^{\epsilon}.
\label{HJ36}
\end{equation}
\label{LHJ8}
\end{Lemma}
{\bf Proof}\\
It's enough to work with $z(x) = x.Du^{\epsilon}(x) =x_i u^{\epsilon}_{x_i}$. Firstly, \vspace{.05in}\\\hh
$z_{x_k}=u^{\epsilon}_{x_k} + x_i u^{\epsilon}_{x_i x_k}$, \vspace{.05in}\\\hh
$z_{x_k x_k}=u^{\epsilon}_{x_k x_k} + u^{\epsilon}_{x_k x_k}+x_i u^{\epsilon}_{x_k x_k x_i}$. \vspace{.05in}\\
Therefore, \vspace{.05in}\\\hh
$Dz=Du^{\epsilon} + x_i Du^{\epsilon}_{x_i}$, \vspace{.05in}\\\hh
$\Delta z=2 \Delta u^{\epsilon}+x_i \Delta u^{\epsilon}_{x_i}$. \vspace{.05in}\\
Besides, differentiate (\ref{HJ21}) with respect to $x_i$
\begin{equation}
DH(Du^{\epsilon}).Du^{\epsilon}_{x_i} =\epsilon \Delta u^{\epsilon}_{x_i}.
\label{HJ37}
\end{equation}
Hence:
\begin{align}
DH(Du^{\epsilon}).Dz - \epsilon  \Delta z &=DH(Du^{\epsilon}).Du^{\epsilon} - 2 \epsilon \Delta u^{\epsilon}+x_i( DH(Du^{\epsilon}).Du^{\epsilon}_{x_i} -\epsilon \Delta u^{\epsilon}_{x_i})
\notag\\
&=DH(Du^{\epsilon}).Du^{\epsilon} - 2 \epsilon \Delta u^{\epsilon}. \notag
\end{align}
We get the lemma.\\\\
This lemma gives us the key idea to find the supersolution of ($\ref{HJ32}$) of the type $z$, and then we can get the result by using Maximum principle.\\
We can choose appropriate $\alpha$, $\beta$ such that $\alpha+\beta >0$ and $\dfrac{2 \alpha+\beta}{\alpha+\beta} = \gamma$.  By using this relation and (H3)
\begin{align}
\label{HJ38}
&DH(Du^{\epsilon}).Dz - \epsilon \Delta z = (\alpha+\beta) (DH(Du^{\epsilon}).Du^{\epsilon} - \gamma \epsilon \Delta u^{\epsilon}) \\
\ge& (\alpha+\beta) (\gamma H(Du^{\epsilon})+\delta - \gamma \epsilon \Delta u^{\epsilon}) =(\alpha+\beta) \delta >0. \qquad \qquad \notag
\end{align}
Let $k = \dfrac{1}{(\alpha+\beta) \delta}$ and let $y(x) = k z(x)+M$ with $M>0$ large enough so that $y|_{\partial U} \ge 0$. Then by (\ref{HJ38}), $y$ is the supersolution of ($\ref{HJ32}$), i.e. 
\begin{equation}
DH(Du^{\epsilon}).Dy - \epsilon  \Delta y \ge 1.
\label{HJ39}
\end{equation}
By the Maximum principle, we easily get:
\begin{equation}
0 \le v^{\epsilon} \le y.
\label{HJ40}
\end{equation}
Therefore, there exists $C>0$ such that $0 \le v^{\epsilon} \le C$.\\
Notice that the boundedness of $U$ plays the crucial role here since it implies the boundedness of $z(x) = \alpha x.Du^{\epsilon}(x) + \beta u^{\epsilon}(x)$. If $U$ is not bounded then $z$ may not be bounded.\\
Like the above section, in order to prove Theorem \ref{THJ3}, we will prove the following theorem
\begin{Theorem}
\label{addTHJ3}
There exists $C>0$ such that
$$
|u^\epsilon_\epsilon(x)| \le C \epsilon^{-1/2}.
$$
\end{Theorem}
{\bf Proof}\\
Differentiate ($\ref{HJ21}$) with respect to $\epsilon$
\begin{equation}
\left\{ \begin{aligned}
DH(Du^{\epsilon}).Du^{\epsilon}_{\epsilon}&= \epsilon \Delta u^{\epsilon}_{\epsilon}+\Delta u^{\epsilon} \h &\mbox{in}&~ U, \vspace{.05in}\\
u^{\epsilon}_{\epsilon}&=0 \h &\mbox{on}&~ \partial U.\\
\end{aligned} \right. 
\label{HJ41}
\end{equation}
There exists $x_0 \in U$ such that
$$
|u^\epsilon_\epsilon(x_0)| = \max_U |u^\epsilon_\epsilon(x)| \ge 0.
$$
Multiply (\ref{HJ41}) by $\sigma^\epsilon$ and then integrate by parts over $U$ as above, we will finally get
$$
u^\epsilon_\epsilon(x_0) = \int_U \Delta u^\epsilon \sigma^\epsilon dx
$$
By Holder's inequality,
$$
|u^\epsilon_\epsilon(x_0) | \le \{ \int_U |D^2 u^\epsilon|^2 \sigma^\epsilon dx\}^{1/2} \{ \int_U \sigma^\epsilon dx\}^{1/2} \le C \epsilon^{-1/2}.
$$
We get the theorem.
\section{Homogenization - The speed of convergence to the effective Hamiltonian}
\h In this section, we study Homogenization and the effective Hamiltonian. We point out the different way of approximation of the Hamilton-Jacobi equation, which includes the vicosity term, and then we study the speed of convergence of the solution of such approximated equation to the effective Hamiltonian.\\\h
Lions, Papanicolaou and Varadhan in \cite {LPV} show the way to find the effective Hamiltonian by considering the following Hamilton-Jacobi equation
\begin{equation}
\epsilon v^{\epsilon} + H(P + Dv^{\epsilon},y) = 0.
\label{HJ42}
\end{equation}
They prove that ${\epsilon}v^{\epsilon}$ converges to $-\bar H(P)$; $v^\epsilon - \min_{\mathbb T^n} v^\epsilon$ converges uniformly to $v$ as ${\epsilon}$ tends to $0$ and $v$ is the viscosity solution of the following cell problem
\begin{equation}
H(P + Dv,y) = \bar H(P),
\label{HJ43}
\end{equation}
where $H\in C^\infty(\mathbb R^n \times \mathbb T^n)$ and satisfies the following conditions 
\begin{equation}
\left\{ \begin{aligned}
|D_xH(p,x)| \le C(1+|p|),~x \in \mathbb T^n,\vspace{.05in}\\
\lim_{|p| \to \infty} \dfrac{H(p,x)}{|p|}= \infty ~\mbox{uniformly in} ~x \in \mathbb T^n.\\
\end{aligned} \right. 
\notag
\end{equation}
Recently, Capuzzo-Dolcetta and Ishii prove the speed of convergence of this problem is $O(\epsilon)$ in \cite{DI}. More precisely, they show that
\begin{equation}
|{\epsilon}v^{\epsilon}+ \bar H(P)| \le C(1+|P|){\epsilon},
\label{HJ43}
\end{equation}
for some constant $C>0$.\\
The proof of the above estimate is really simple and only based on some comparison principles. However, there are still some difficult issues remaining. The most difficult one is that even though $\bar H(P)$ is unique, it's hard to study the encoded information in $\bar H$, especially in the context of weak KAM theory. Also in practice, it's hard to calculate the solution of (\ref{HJ42}).\\
Our approach here is different. Firstly, let's look at the regularized equation with viscosity term of (\ref{HJ42})
\begin{equation}
{\epsilon}v^{\epsilon, \delta}+ H(P+Dv^{\epsilon, \delta},y)= \delta  \Delta v^{\epsilon, \delta}.
\label{HJ44}
\end{equation}
As we have already proved, there exists a constant $C>0$ independent of $\epsilon, \delta$ such that
\begin{equation}
|{\epsilon}v^{\epsilon, \delta}-{\epsilon}v^{\epsilon}| \le C {\delta}^{1/2}.
\label{HJ45}
\end{equation}
Combining (\ref{HJ43}) and (\ref{HJ45}),
\begin{equation}
|{\epsilon}v^{\epsilon, \delta}+ \bar H(P)| \le C {\delta}^{1/2} + C  \epsilon,
\label{HJ46}
\end{equation}
In particular, if we choose $\delta = \epsilon^2$ then
\begin{equation}
|{\epsilon}v^{\epsilon, \epsilon^2}+ \bar H(P)| \le  C  \epsilon.
\label{HJ47}
\end{equation}
This is the motivation for us to consider a slightly different regularized problem as following\\
 Let $z^{\theta}$ be the solution of the following Hamilton-Jacobi equation
\begin{equation}
{\theta}z^{\theta}+ H(P+Dz^{\theta},y)= \theta^2  \Delta z^{\theta}.
\label{HJ48}
\end{equation}
We will show in this section that
\begin{Theorem}
\label{effective}
There exists a constant $C>0$ such that
$$
||\theta z^\theta +\bar H (P)||_{L^\infty} \le C \theta.
$$
\end{Theorem}
Similar to the the previous sections, we will show
\begin{equation}
|\dfrac{\partial (\theta z^{\theta})}{\partial \theta}| \le C.\notag
\end{equation}
Gomes in \cite{G} also considers some equations similar to (\ref{HJ48}) to study the properties of Mather measures.\\
Although our method is slightly complicated than the method in \cite{DI}, it creates a constructive way to study the effective Hamiltonian and we will use it to study weak KAM theory and Mather measures elsewhere in the future.\\\\
Now we will study the properties of $z^\theta$ and prove Theorem \ref{effective}.\\ 
Firstly, we have some following standard observations: $z^{\theta}$ is unique hence $\mathbb T^n$-periodic and from \cite{L}, there exists a constant $C>0$ such that
\begin{equation}
||{\theta}z^{\theta}||_{L^{\infty}},~||Dz^{\theta}||_{L^{\infty}} \le C. \notag
\end{equation}
Also from the $\mathbb T^n$-periodic property of $z^\theta$ and the boundedness of $||Dz^{\theta}||_{L^{\infty}}$, we have furthermore that
\begin{equation}
|z^{\theta}(x)-z^{\theta}(y)| \le C ||Dz^{\theta}||_{L^{\infty}} \le C \h \forall~ x,y \in \mathbb R^n.\notag
\end{equation}

\begin{Lemma} Let $w^{\theta}= \dfrac{|Dz^{\theta}|^2}{2}$ then
\begin{equation}
2{\theta}w^{\theta}+ D_pH.Dw^{\theta}+D_xH.Dz^{\theta}= \theta^2  \Delta w^{\theta}-\theta^2  |D^2z^{\theta}|^2.
\label{HJ49}
\end{equation}
\label{LHJ9}
\end{Lemma}
The equation here is of first type since ($\ref{HJ49}$) contains $w^{\theta}$. Using the same method as in Section 2, we introduce the fake time-dependent function $v^{\theta}$ such that $v^{\theta}(x,t)=e^t w^{\theta}(x)$ for $t \in [0,T]$ for some $T>0$ fixed.\\
 Then $v^{\theta}$ satisfies
\begin{equation}
2{\theta}v^{\theta}_t+ D_pH.Dv^{\theta}+e^t D_xH.Dz^{\theta}= \theta^2  \Delta v^{\theta}-\theta^2 e^t  |D^2z^{\theta}|^2.
\label{HJ50}
\end{equation}
{\bf Adjoint method:} We now introduce the adjoint equation of ($\ref{HJ50}$):
\begin{equation}
\left\{ \begin{aligned}
-2\theta \sigma^{\theta}_t - \mbox{div}(D_pH\sigma^{\theta})&= \theta^2 \Delta \sigma^{\theta} \h (x,t) \in \mathbb R^n \times (0,T)\vspace{.05in}\\
 \sigma^{\theta}|_T&=\delta_{x_0}.\\
\end{aligned} \right. 
\label{HJ51}
\end{equation}
By the theory of distributions (see Chapter 5 in \cite{RR}), $\sigma^\theta$ is unique and $\sigma^\theta \in C^\infty(\mathbb R^n \times (0,T))$.\\
Similar to section 2 above, we have some properties of $\sigma^{\theta}$ as following
\begin{Lemma} {\bf Properties of $\sigma^{\theta}$}
\label{LHJ10}
\begin{itemize}
\item[(i)] $\sigma^{\theta}(x,t) \ge 0$ ~for $ (x,t) \in \mathbb R^n \times (0,T)$,
\item[(ii)] $\int_{\mathbb R^n} \sigma^{\theta}(x,t)dx=1$ ~for $ t \in (0,T)$.
\end{itemize}
\end{Lemma}

\begin{Lemma}
There exists a constant $C>0$ such that
\begin{equation}
\theta^2 \int_0^T \int_{\mathbb R^n} |D^2 z^{\theta}|^2 \sigma^{\theta} \le C.
\label{HJ52}
\end{equation}
\label{LHJ11}
\end{Lemma}
Again, all the estimates here don't depend on the choice of $x_0$ as stated carefully in section 2.
\begin{Theorem}
There exists a constant $C>0$ such that
\begin{equation}
|(\theta z^{\theta})_{\theta}(x)| \le C.
\label{HJ53}
\end{equation}
\label{THJ4}
\end{Theorem}
{\bf Proof}\\
Firstly, differentiate ($\ref{HJ48}$) with respect to $\theta$ 
\begin{equation}
z^{\theta}+ \theta z^{\theta}_{\theta}+D_pH.Dz^{\theta}_{\theta}=\theta^2  \Delta z^{\theta}_{\theta} + 2 \theta \Delta z^{\theta}.
\label{HJ54}
\end{equation}
Doing the same steps as in Theorem $\ref{THJ2}$, we get the following
\begin{align}
\label{HJ55}
2{\theta}(e^T z^{\theta}_{\theta}(x_0)-\int_{\mathbb R^n} z^{\theta}_{\theta}(x) \sigma^{\theta}(x,0)dx)- \theta \int_0^T\int_{\mathbb R^n} z^{\theta}_{\theta} \sigma^{\theta}dx dt ~+ \qquad \\
+\int_0^T\int_{\mathbb R^n} e^t z^{\theta} \sigma^{\theta}dx dt=2\theta \int_0^T\int_{\mathbb R^n} e^t \Delta z^{\theta} \sigma^{\theta}dx dt .\notag
\end{align}
Let
\begin{align}
A&= 2{\theta}(e^T z^{\theta}_{\theta}(x_0)-\int_{\mathbb R^n} z^{\theta}_{\theta}(x) \sigma^{\theta}(x,0)dx)- \theta \int_0^T\int_{\mathbb R^n} z^{\theta}_{\theta} \sigma^{\theta}dx dt,\quad \notag \\
B&= \int_0^T\int_{\mathbb R^n} e^t z^{\theta} \sigma^{\theta}dx dt. \qquad \qquad\qquad\qquad\qquad\qquad\qquad\qquad\quad\notag
\end{align}
Notice that $|A+B| \le C$ for some positive constant $C$ independent of the choice of $x_0$ by Lemma $\ref{LHJ11}$.\\ 
We have two following observations
\begin{itemize}
\item[(i)] Take any $x' \in \mathbb R^n$, we can control $B$ in term of $z^{\theta}(x')$ by using the property $|z^{\theta}(x) - z^{\theta} (x')| \le C$ for all $x \in \mathbb R^n$. More explicitly,
\begin{equation}
|B- (e^T-1) z^{\theta} (x')| \le  \int_0^T\int_{\mathbb R^n} e^t |z^{\theta}(x)-z^{\theta}(x')| \sigma^{\theta}dx dt \le C(e^T-1)=C.
\label{HJ56}
\end{equation}
\item[(ii)] By using the same arguments as in section 2, we may assume without loss of generality that there exist $x_1, x_2 \in \mathbb R^n$ such that
\begin{equation}
z^{\theta}_{\theta}(x_1) =m=  \min_{\mathbb R^n}z^{\theta}_{\theta}(x) \le z^{\theta}_{\theta}(x) \le \max_{\mathbb R^n}z^{\theta}_{\theta}(x) = M = z^{\theta}_{\theta}(x_2).
\notag
\end{equation}
\end{itemize}
Let $\sigma^{\theta,x_1}$ be the solution of the adjoint equation corresponding to $x_1$ (i.e. we let $x_0=x_1$ in the adjoint equation (\ref{HJ51})), then
\begin{equation}
A= 2{\theta}(e^T z^{\theta}_{\theta}(x_1)-\int_{\mathbb R^n} z^{\theta}_{\theta}(x) \sigma^{\theta,x_1}(x,0)dx)- \theta \int_0^T\int_{\mathbb R^n} z^{\theta}_{\theta} \sigma^{\theta,x_1}dx dt \le (e^T-1) \theta z^{\theta}_{\theta}(x_1).
\notag
\end{equation}
Therefore, by both of the observations above,
\begin{equation}
-C \le A+B \le (e^T-1) \theta z^{\theta}_{\theta}(x_1) + B \le (e^T-1) m \theta + (e^T-1) z^{\theta} (x') +C,
\notag
\end{equation}
which implies that
\begin{equation}
-C \le (e^T-1) (m \theta +  z^{\theta} (x') ) \quad \forall ~x' \in \mathbb R^n.
\label{HJ59}
\end{equation}
Similarly, let $\sigma^{\theta,x_2}$ be the solution of the adjoint equation corresponding to $x_2$ (i.e. we let $x_0=x_2$ in the adjoint equation (\ref{HJ51})), then
\begin{equation}
A= 2{\theta}(e^T z^{\theta}_{\theta}(x_2)-\int_{\mathbb R^n} z^{\theta}_{\theta}(x) \sigma^{\theta,x_2}(x,0)dx)- \theta \int_0^T\int_{\mathbb R^n} z^{\theta}_{\theta} \sigma^{\theta,x_2}dx dt \ge (e^T-1) \theta z^{\theta}_{\theta}(x_2).
\notag
\end{equation}
Therefore, the following estimate holds
\begin{equation}
(e^T-1) (M \theta +  z^{\theta} (x') ) \le C \quad \forall ~x' \in \mathbb R^n.
\label{HJ61}
\end{equation}
Combining ($\ref{HJ59}$) and ($\ref{HJ61}$)
\begin{equation}
(e^T-1) |\theta z^{\theta}_{\theta} (x)+z^{\theta} (x)| \le C \quad \forall ~x \in \mathbb R^n.
\label{HJ62}
\end{equation}
We get the theorem.
\section{Appendix}
We will prove the uniqueness of the solution $u^{\epsilon}$ of equation ($\ref{HJ21}$).
\begin{Theorem}
If $u$ and $v$ are the solutions of ($\ref{HJ21}$) then we get $u=v$.
\label{AHJ1}
\end{Theorem}
{\bf Proof}\\
It's enought to prove that $ u \le v$. \vspace{.05in}\\
If we have $H(Du) - \epsilon \Delta u < H(Dv) - \epsilon \Delta v$ in $U$ and $ u \le v$ on $\partial U$ then we easily get $u \le v$ in $U$ by the usual Maximum principle.\\
The strategy is to find a sequence of functions $\{z^{\theta}\}$ such that $z^{\theta}$ converges uniformly to $v$ and 
$$
H(Du) - \epsilon \Delta u < H(Dz^{\theta}) - \epsilon \Delta z^{\theta}~\mbox{in}~U;~\mbox{and}~  u \le z^{\theta}~\mbox{on}~ \partial U.
$$
Hence we will get $u \le z^{\theta}$ for all $\theta$, which implies $u \le v$.\\
By Remark $\ref{RHJ1}$, for $t>1$ we have
\begin{equation}
H(tp) \ge t^{\gamma} H(p) + \dfrac{\delta}{\gamma+1}(t^{\gamma}-1).
\notag
\end{equation}
Let $z= sv+t(x.Dv+M)$ where $M>0$ is to be a suitable chosen constant.\\ 
We can see that the function $z$ here is similar to the one in Lemma $\ref{LHJ8}$. We have: \vspace{.05in}\\\h
$Dz = (s+t)Dv + tx_i Dv_{x_i}$, \vspace{.05in}\\\h
$\Delta z = (s+2t)\Delta v + t x_i \Delta v_{x_i}$. \vspace{.05in}\\
For $s$ close to $1$, for $t>0$ close to $0$ and $s+t>1$,
\begin{align}
&H(Dz) - \epsilon \Delta z = H((s+t)Dv + tx_i Dv_{x_i}) - \epsilon (s+2t)\Delta v + t \epsilon x_i \Delta v_{x_i}  \vspace{.05in}\qquad\qquad\qquad\notag\\
=& H((s+t)Dv) +tDH((s+t)Dv).(x_i Dv_{x_i})+t^2O(1) -  \epsilon (s+2t)\Delta v + t \epsilon x_i \Delta v_{x_i}  \vspace{.05in}\quad\notag\\
=&H((s+t)Dv) +tDH(Dv).(x_i Dv_{x_i})+t((s+t)-1)O(1)+t^2O(1)- \qquad\qquad\qquad\notag \vspace{.05in} \\
&-  \epsilon (s+2t)\Delta v + t \epsilon x_i \Delta v_{x_i} \qquad\notag \vspace{.05in}\\
\ge& (s+t)^{\gamma}H(Dv)+\dfrac{\delta}{\gamma+1}((s+t)^{\gamma}-1)+t((s+t)-1)O(1)+t^2O(1)-\qquad\qquad\qquad\notag \vspace{.05in}\\
&-  \epsilon (s+2t)\Delta v.\qquad  \vspace{.05in}
\notag
\end{align} 
For $\theta>0$, let $t=(1+\theta)^{\gamma} - (1+\theta)$ and $s=2(1+\theta)-(1+\theta)^{\gamma}$.\\ 
Let $z^\theta= sv+t(x.Dv+M)$ corresponding to $s,t$ chosen above.\\ 
Notice that $z^\theta$ converges uniformly to $v$ as $\theta$ tends to $0$.\\
Furthermore, $(s+t)^{\gamma}=s+2t=(1+\theta)^{\gamma}$ and for $\theta$ small enough
\begin{equation}
\dfrac{\delta}{\gamma+1}((s+t)^{\gamma}-1)+t((s+t)-1)O(1)+t^2O(1) >0.
\notag
\end{equation}
Hence we get $H(Dz^\theta) -\epsilon \Delta z^\theta >0$.\\
Finally, choose $M$ large enough to guarantee $z^\theta \ge u$ on $\partial U$, we get the theorem.
\begin{acknowledgement}
I would like to express my appreciation to my advisor, Lawrence C. Evans for giving me the problems and plenty of fruitful discussions. I thank Scott Armstrong, Filippo Cagnetti, Charlie Smart for their helpful discussions and suggestions. Finally, I would like to thank the anonymous referee for his kind comments and suggestions.
\end{acknowledgement}

\end{document}